\documentclass[12pt]{amsart}
\usepackage{subfigure}
\usepackage{graphicx}
\usepackage{rotating}
\usepackage{morefloats}
\usepackage[utf8]{inputenc}
\usepackage[T1]{fontenc}
\usepackage{amssymb}
\usepackage{amsmath}    
\def\R{I\kern -0,37 em R}
\def\P{I\kern -0,37 em P}
\def\Z{I\kern -0,37 em Z}

\begin{document}
\title[local and global equivalence]{On the local equivalence of partial differential equations}
\author{A. Kumpera}
\address[Antonio Kumpera]{Campinas State University, Campinas, SP, Brazil}
\email{antoniokumpera@hotmail.com}

\date{August, 2016}


\maketitle

\begin{abstract}
In all the practical applications of partial differential equations, what is mostly needed and what is in fact hardest to obtains are the solutions of the system or, occasionally, some specific solutions. This work is based on a most enlightening \textit{Mémoire} written by Élie Cartan in 1914 and that the majority ignores (\cite{Cartan1914}, \textit{see also} \cite{Cartan1937}). We discuss a setting for the local equivalence problem and illustrate it by some examples. It should also be noted that any integration process or method is in fact a local equivalence problem involving a suitable model.
\end{abstract}

\section{Introduction}
We discuss here the equivalence problem in the realm of Ehresmann's prolongation spaces of which, jet spaces are a special case. All the emphasis is given to the first order contact system canonically associated to a given partial differential equation since the local or, inasmuch, the global equivalence of these Pfaffian systems will entail the corresponding equivalence for the equations. A fundamental ingredient in our approach is the consideration of \textit{merihedric} prolongation spaces (\textit{prolongements mériédriques}) as defined by Élie Cartan in many of his writings (\textit{cf}. the references). Though rather absent in the recent literature, the reader will find interesting examples in \cite{Kumpera1999},  \cite{Kumpera2014} and \cite{Kumpera2015}. As for the calculations, we try to reduce them to the strict minimum. The reason for considering the local equivalence problem relies on the eventual possibility of exhibiting a local model equivalent to the given equation and ultimately hope that this model is simple enough so as to become  integrable by some known methods. In this respect, some information is provided in \cite{Cartan1910} and \cite{Kumpera2016}. As for most (perhaps all) equivalence problems, we seek to determine \textit{a fundamental set of invariants} associated to a differential system that will characterise as well as describe the local equivalences.

\section{Generalities}
In this section we recall, briefly, some very standard facts so as to fix the notations and subsequently discuss, in more detail, our present aims and techniques.

\vspace{5 mm}
\noindent
As is usual, a \textit{fibration} is a surjective map $P\longrightarrow M$ of maximum rank and we define a \textit{multi-fibration} as being a sequence, finite or infinite in length,

\begin{equation*}
\cdots~\longrightarrow~P_\mu~\longrightarrow~P_{\mu-1}~\longrightarrow~\cdots~\longrightarrow~P_1\longrightarrow P_0 ~,
\end{equation*}

\vspace{3 mm}
\noindent
where all the projections $\rho_{\beta,\alpha}:P_\alpha\longrightarrow P_\beta$ , $\alpha\geq\beta$, are fibrations, the map $\rho_{\beta,\alpha}$ being the composite of the successive projections. Following Cartan, we say that a multi-fibration is a \textit{merihedric} prolongation space with respect to the basic fibration $P_1\longrightarrow P_0$ when there exists a prolongation algorithm $\wp$ such that any local transformation $\varphi:U\longrightarrow V$ on the manifold $P_1$ that preserves the transversality with respect to the fibres of $\pi:P_1\longrightarrow P_0$ can be prolonged (extended) to a local transformation $\varphi_{\mu}=\wp_{\mu}(\varphi)$ operating on the manifold $P_{\mu},$ the usual commutativity relations with respect to the projections as well as to the composition of transformations and the passage to the inverses being preserved. However, we shall only require, with the obvious notations, that the open sets $U_{\mu}$ and  $V_{\mu}$ project surjectively onto $U_{\nu}$ and $V_{\nu}$ whenever $\nu\leq\mu$ without requiring that the former be the inverse images of the latter. The preservation of transversality simply means that any linear sub-space transversal to the fibres of $\pi$ is transformed by $\varphi_*$ onto a transversal sub-space at the target point. We next require that to every local section $\sigma$ of $\pi$ corresponds, at each level $\mu,$ a \textit{prolonged} section $\wp_{\mu}\sigma$ and that these sections project one upon the other. Finally we require that, at each level $P_\mu,~\mu\geq 2,$ a Pfaffian system $\mathcal{C}_{\mu}$ be defined and that it satisfies the following properties:

\vspace{4 mm}
(a) $\rho^*_{\mu,\mu-1}\mathcal{C}_{\mu-1}\subset\mathcal{C}_{\mu},$

\vspace{3 mm}
(b) for any local section $\sigma$ of $\pi,$ the image of $\wp_{\mu}\sigma$ is an integral sub-manifold of $\mathcal{C}_{\mu}$ and, conversely, any section $\tau$ of the \textit{source} projection (the iterated composition of all the projections starting with $P_\mu~\longrightarrow~P_{\mu-1}$) whose image is an integral of $\mathcal{C}_{\mu}$ is of the form $\wp_{\mu}\sigma$ and 

\vspace{3 mm}
(c) the $\mu-$th prolongation of any local transformation of $P_1$ is an automorphism of $\mathcal{C}_{\mu},$ the converse being verified just locally.

\vspace{4 mm}
\noindent
\textbf{Remark.} It would be much more convenient to work on $k-$th order Grassmannian bundles rather than on $k-$th order Jet spaces since this would avoid transversality concerns. Nevertheless, Ehresmann's Jet bundles are very convenient for calculations that frequently cannot be avoided.

\vspace{4 mm}
\noindent
A similar statement holds for vector fields (\textit{infinitesimal transformations}) and their prolongations where the operations are now the linear operations together with the Lie derivative and the bracket of vector fields. We also require that this \textit{infinitesimal}  prolongation procedure be compatible with the \textit{finite} prolongation of the elements of a local $1-$parameter group. Since such a $1-$parameter group always initiates with the Identity transformation, the above transversality assumption is entirely irrelevant in the infinitesimal context. 

\vspace{5 mm}
\noindent
For convenience, we denote by $\pi:P\longrightarrow M$ the bottommost fibration $P_1~\longrightarrow~P_0,$ call $P$ the \textit{total space} and $M$ the \textit{base space}. The most requested prolongation spaces are, of course, the sequences of Ehresmann's jet spaces of all orders associated to given fibrations $\pi,$ where we consider $k-$jets of local sections. We shall not provide any details on jet spaces since these are a well known subject. Nevertheless, we call the attention of the reader that, contrary to the usual procedure, the index $k,$ or its equivalent, will always be placed as a sub-script, that the $k-$jet of a local section $\sigma,$ at the point $x,$ will be written $j_k\sigma(x)$ and that the $k-$th order prolongation procedure will be denoted by $p_k(~).$

\vspace{5 mm}
\noindent
A \textit{partial differential equation} of order $k$ is a sub-fibration $\beta_k:\mathcal{S}~\longrightarrow~P$ of the target projection  $\beta_k:J_k\pi~\longrightarrow~P,$ where $J_k\pi$ denotes the space of all the $k-$jets of local sections of the fibration  $\pi.$ Occasionally, such a differential equation can just be defined over an open subset of $P.$ A \textit{solution} of $\mathcal{S}$ is, of course, a local section $\sigma$ of $\pi$ such that its $k-$jet $j_k\sigma$ becomes a section of the \textit{source} projection $\alpha_k:\mathcal{S}~\longrightarrow~M.$ A solution of $\mathcal{S}$ is also a solution of any prolongation of $\mathcal{S}.$ The $k-$th order \textit{contact} system $\mathcal{C}_k$ is the Pfaffian system defined on $J_k\pi$ and generated by all the Pfaffian $1-$forms annihilating the tangent spaces to any \textit{holonomic} section \textit{i.e.}, a section of the form $j_k\sigma.$ Two fundamental properties of this Pfaffian system state that (a) any section $\tau$ of the source projection $\alpha_k:J_k\pi~\longrightarrow~M$ whose image is an integral of $\mathcal{C}_k$ is holonomic and (b) the prolongation of any local transformation of $P$ is an automorphism of $\mathcal{C}_k$ though the converse just holds locally. This family of Pfaffian systems plays therefore the role of the Pfaffian systems defined on the spaces belonging to a merihedric prolongation space. The $k-$th order \textit{contact} system $\mathcal{C}_k\mathcal{S}$ associated to the equation $\mathcal{S}$ is the restriction of $\mathcal{C}_k$ to the sub-manifold $\mathcal{S}.$

\vspace{5 mm}
\noindent
For the Ehresmann prolongation space built up in terms of the jet spaces, we have the additional bonus of a prolongation procedure for the differential equations. We can, however, also devise such a prolongation procedure for any merihedric prolongation space by adopting exactly the same definitions. In fact, let us take a sub-manifold $\mathcal{M}\subset P_{\mu}$ and let us define the $\ell-$th prolongation of $\mathcal{M}$ as the sub-set of $P_{\mu+\ell}$ composed by all the elements $\wp_{\mu+\ell}\sigma(x)$ such that $\wp_{\mu}\sigma(x)\in\mathcal{M}$ and, furthermore, such that the image of the section $\wp_{\mu}\sigma$ is tangent up to order $\ell$ to the sub-manifold $\mathcal{M}$ at the point $\wp_{\mu}\sigma(x).$ All the properties of the standard prolongation algorithm for differential equations transcribe in this more general context and, in particular, the $\ell-$prolongation of $\mathcal{M}$ can be obtained inasmuch by iterated first order prolongations. Unfortunately, the prolongation of $\mathcal{M}$ is not necessarily, as in the case of differential equations, a differentiable sub-manifold of $P_{\mu+\ell}.$

\vspace{5 mm}
\noindent
In studying the equivalence, local or global, of two differential equations we shall always consider both equations defined over the same basic fibration $\pi$ for, if the two equations were situated above distinct basic fibrations, we could always shift one of them onto the basic fibration of the other by simply taking a local or global diffeomorphism among the two fibrations (a diffeomorphism respecting the fiberings) and place all the data on one of them. Concerning the basic properties and operations relative to differential equations and their prolongations, we refer the reader to \cite{Kumpera1972}. 

\vspace{5 mm}
\noindent
Let us now begin by examining the notion of \textit{equivalence} for differential equations. This notion must respect the solutions of the equations \textit{i.e.}, a solution of one of them must be transformed into a solution of the other hence must respect the sophisticated structure of jet spaces.

\vspace{2 mm}
\newtheorem{lagrange}[DefinitionCounter]{Definition}
\begin{lagrange}
Two $k-$th order differential equations $\mathcal{S}$ and $\mathcal{S}'$ are said to be locally absolutely equivalent in the neighborhoods of $X\in\mathcal{S}$ and $X'\in\mathcal{S}'$ when there exists a local transformation $\varphi:U~\longrightarrow~U'$ on the total manifold P such that its prolongation $p_k\varphi$ transforms X in X' and, further, becomes a diffeomorphism when restricted to the appropriate open subsets of $\mathcal{S}$ and $\mathcal{S}'.$ The equivalence is global when $\varphi$ is global. 
\end{lagrange}

\vspace{2 mm}
\noindent
Obviously, $p_k\varphi^{-1}$ is also an equivalence and any one of them transforms solutions into solutions. If, for instance, $\sigma$ is a solution of $\mathcal{S},$ then $\varphi\circ\sigma\circ\overline{\varphi}^{-1}$ is a solution of $\mathcal{S}',$ where $\overline{\varphi}=\pi\circ\varphi\circ\sigma.$  

\vspace{5 mm}
\noindent
The notion of absolute equivalence seems too much restricted when confronted with applications and we now give, below, a broader and more suitable definition. Needless to say that both definitions are stated in Cartan's \textit{Mémoire} \cite{Cartan1914}.

\vspace{2 mm}
\newtheorem{laguère}[DefinitionCounter]{Definition}
\begin{laguère}
Two differential equations are said to be locally equivalent when they admit prolongations of a certain order that are locally absolutely equivalent. Global equivalence has a similar definition.
\end{laguère}

\vspace{5 mm}
\noindent
Most probably, the reader already observed that the last definition is a big cheat. In fact, it is so since the absolute equivalence of some prolongations implies, by projection, the absolute equivalence of the given equations. Nevertheless, quite often it is easier and more convenient to deal with prolongations than with the given equations. Unfortunately, Cartan is not here to give us a more convincing argument.

\vspace{5 mm}
\noindent
The first example occurs in the case of \textit{ordinary differential equations}. In this case, when the image spaces have the same dimension (the same number of dependent variables), any two equations are always locally equivalent at non-singular points. When singularities are present, the matter becomes considerably more involved. As a second example, we mention that all the integrable Pfaffian systems with the same rank and co-rank are locally equivalent. As for the global equivalence, in the present case as well as in the case of ordinary differential equations, this is a beautiful and still unresolved problem in topology. A very interesting context can be found in \cite{Almeida1985}. Let us now have a glance at non-integrable Pfaffian systems, a rather nice example being provided by the \textit{Flag Systems}. These are systems of co-rank equal to 2 (when the characteristics vanish) hence the integral manifolds are just $1-$dimensional curves, \textit{l'intégrale générale ne dépendant que de fonctions arbitraires d'une seule variable}, as would claim Cartan. However, the reader might delight himself in examining the many local models exhibited in \cite{Kumpera2014}.

\vspace{5 mm}
\noindent
Let us now go one step further and define the \textit{merihedric equivalence}. This equivalence - actually, there are uncountably many possibilities - differs from the above \textit{holonomic} equivalence but not too much. As Cartan himself did say \textit{nous allons voir que c'est au fond le seul} (the holonomic prolongation, p. 6, $l.$ 15). What we shall do is to adapt, to the jet spaces sequence, a merihedric process namely, we consider the following prolongation space where $k$ is the order of the equations:

\begin{equation*}
\cdots~\longrightarrow~P_{k+\ell}~\longrightarrow~\cdots~\longrightarrow~P_{k+1}~\longrightarrow~J_k\pi~\longrightarrow~\cdots~\longrightarrow~P~\longrightarrow~M ~.    
\end{equation*}

\vspace{3 mm}
\noindent
The definitions of equivalence are exactly the same as above since we also have available a prolongation process, for any differential equation (a sub-manifold) contained in $J_k\pi,$ in the merihedric case.

\vspace{5 mm}
\noindent
We now return to the special case of Flag Systems that provide a very nice illustration of the merihedric situation and, further, are prolongation spaces of finite length. Before doing so, we recall that in \cite{Kumpera1999}, we mention that the standard \textit{total derivative}, in jet spaces, can be trivially extended to a \textit{total Lie derivative} for differential forms and, with this in mind, we shall give attention to Pfaffian systems instead of partial differential equations, which was by far Cartan's preference. In this respect, we also recall that the total Lie derivative of the $k-$th order canonical contact system, defined on $J_k\pi,$ is precisely the corresponding contact system on the $(k+1)-$st level. Let us then return to the sequence

\begin{equation*}
(P_\ell,S)~\longrightarrow~(P_{\ell-1},\overline{S}_1)~\longrightarrow~\cdots~\longrightarrow~(P_0,\overline{S}_{\ell}) ~,
\end{equation*}

\vspace{3 mm}
\noindent
exhibited in \cite{Kumpera2014}, p. 8, $l.$ -3, where \textit{S} is a Flag System of length $\ell$ defined on the space $P=P_{\ell},$ with null characteristics and where the $\overline{S}_{\ell-\nu}$ denote Pfaffian systems canonically isomorphic to the successive derived systems of \textit{S} though defined on more appropriate spaces so as to have also null characteristics. The above finite sequence is a prolongation space together with Pfaffian systems, at each level, naturally associated to the initially given system \textit{S}. The terminating system $\overline{S}_\ell$ vanishes and $\overline{S}_{\ell-1}$ is a  Darboux system in $3-$space. As shown in \cite{Kumpera2014}, every local or infinitesimal automorphism of this Darboux system extends (prolongs) canonically to an automorphism of $(P_{\ell-\nu},\overline{S}_{\nu}),$ this correspondence becoming, moreover, an isomorphism of pseudo-groups as well as of pseudo-algebras. We are not to be concerned with local sections since $P_1$ can be considered as the total space, whereas the base space $P_0$ collapses to an open set, eventually to a point. We claim that each Pfaffian system $\overline{S}_{\nu}$ is the prolongation of $\overline{S}_{\nu+1}.$ As is, not only it seems that we are walking backwards but worse, the prolongation algorithm as described previously is somehow missing. Nevertheless, in the passage from $\overline{S}_{\nu+1}$ to $\overline{S}_{\nu},$ three options are possible (\cite{Kumperaa1982}) and, upon choosing the appropriate one, we can thereafter apply the prolongation method as described earlier. Furthermore, if we choose the appropriate coordinates so as that the form $\omega=dy-zdx$ generates the Darboux system $\overline{S}_{\ell-1},$ then total derivatives reappear again and Cartan, as always, is perfectly right in his claims.  

\section{The local equivalence of differential systems}
We shall now examine a few aspects of this local equivalence problem where a most enlightening example is provided by the study of the local equivalence of geometrical structures since, in most cases, this is performed by studying the local equivalence of the defining differential equations for these structures (\cite{Kumpera2017}).

\vspace{5 mm}
\noindent
$(a)~Ordinary~Differential~Equations.$

\vspace{3 mm}
\noindent
An ordinary differential equation of order 1 in \textit{n} unknown functions is simply a vector field on a manifold \textit{P}, where we take the fibratiom $\pi:P~\longrightarrow~\mathbf{R}$ with $dim~P=n+1.$ A $k-$th order equation is a contact vector field (infinitesimal automorphism of the corresponding contact structure) defined on $J_k\pi.$ As mentioned earlier, any two ordinary differential equations are always locally equivavalent at two non-singular points. As for the global equivalence, this is a topological problem that involves the nature of the global solutions and, of course, compactness is the main ingredient. We also observe that two ordinary differential equations of distinct orders can be locally, \textit{viz.} globally, equivalent. As for the local equivalence in the neighborhoods of singular points, there is an immense literature on this subject hence we shall not enter here into details concerning this matter.

\vspace{3 mm}
\noindent
$(b)~Linear~Partial~Differential~Equations.$

\vspace{3 mm}
\noindent
Such equations are defined, most conveniently, as vector sub-bundles of the $k-$th order Jet spaces of local sections of vector bundles \textit{i.e.}, the \textit{total space} of the initial fibration is the total space of a vector bundle and, consequently, the corresponding Jet space is as well a vector bundle. In this case, it is most appropriate to restrict the local equivalences to linear morphisms though two linear equations can actually be equivalent via a non-linear local transformation. A first systematic investigation on the local equivalence of linear differential systems was initiated by Drach, Picard and Vessiot, whereupon resulted the well known \textit{Théorie de Picard-Vessiot} (\textit{cf.} \cite{Kumpera2016}).

\vspace{3 mm}
\noindent
$(c)~General~Differential~Equations.$

\vspace{3 mm}
\noindent
Who undoubtedly most contributed not only with results but mainly with ideas and methods, in this general equivalence problem, was Élie Cartan and the reader is invited to examine the table of contents of the \textit{O\!euvres Complètes, Partie II}. Very significant contributions were also brought by Bernard Malgrange, in particular for the case of elliptic equations, as well as by Masatake Kuranishi for the case of the defining equations of a Lie pseudo-group namely, those equations that have the additional algebraic structure of a groupoid. Here again, we shall not enter into further details since the subject is by far too extensive.

\vspace{3 mm}
\noindent
$(d)~Pfaffian~Systems.$

\vspace{3 mm}
\noindent
Initial contributions were, of course, brought by Pfaff himself and later expanded by contemporary Sophus Lie (\textit{Gasammelte Abhandlungen, Vol. 6}) and Mark von Weber. But again, the most striking contributions were placed forward by no other than Élie Cartan who gave us much insight in the local equivalence of the non-integrable systems. As mentioned earlier, any two integrable Pfaffian systems are \textit{everywhere} locally equivalent if and only if they have the same rank and co-rank (the underlying manifolds ought to be of the same dimension). It is also a straightforward consequence of linearity in the tangent and co-tangent spaces that any two Pfaffian systems, integrable or not, are locally equivalent to first order if and only if their ranks and co-ranks are equal.   

\vspace{3 mm}
\noindent
$(e)~Exterior~Differential~Systems.$

\vspace{3 mm}
\noindent
All the merits belong here to Élie Cartan who, in fact, introduced such systems in the study of various geometrical problems (\cite{Cartan1945}). Nevertheless, important contributions are also due to Masatake Kuranishi.

\section{The local equivalence of Lie groupoids}
Since we are interested in the local equivalence of partial differential equations, the sole Lie groupoids to be considered in the sequel are those whose total spaces are sub-manifolds in some jet space and we begin here with some general considerations.

\vspace{5 mm}
\noindent
Given the fibration $\pi=p_1:P=M\times M~\longrightarrow~M$ and the corresponding $k-$th order Jet bundle $J_k\pi,$ we denote by $\Pi_kM$ the groupoid of all the $k-$jets of invertible sections of $\pi$ (\cite{Kumpera2015}). In other terms, $\Pi_kM$ is the set of all the invertible $k-$jets with source and target in \textit{M} and is an open dense subset of $J_k\pi.$ The Lie groupoids that we shall be interested in are the sub-groupoids of $\Pi_kM$ that are inasmuch sub-manifolds (not necessarily regularly embedded) of the ambient groupoid and, together with the differentiable sub-manifold sructure, become Lie groupoids in the sense of Ehresmann (\textit{cf.} the above citation). Such $k-$th order Lie groupoids will, most often, be denoted by $\Gamma_k.$ We next observe that the standard prolongation of order \textit{h} of a Lie groupoid, in the sense of partial differential equations, is again a Lie groupoid at order $k+h$ and all the usual properties relative to differential equations transcribe for the Lie groupoids. Moreover, it can be shown that the above-mentioned properties still hold in the more general context of merihedric prolongations. In considering the local equivalences, we proceed in analogy with what is usually done in Lie group theory and restrict our attention to open neighborhoods of the unit elements of the groupoids. More precisely, we look for conditions under which there exists a local contact transformation that maps isomorphically an open neighborhood of the unit elements in one of the groupoids onto a similar neighborhood of the other. The term "isomorphically" means that it preserves the algebraic as well as the differentiable properties and "an open neighborhood of the unit elements" refers to a neighborhood of an open set in the units sub-manifold. Since a contact transformation, as above, operating on the sub-manifolds is by necessity induced by a local contact transformation operating in the ambient Jet space, we ultimately consider, in our setting, a local contact transformation $\varphi:\mathcal{U}\longrightarrow\mathcal{U}'$ that transforms accordingly the respective intersections with the groupoids and where $\mathcal{U}$ and $\mathcal{U}'$ are open sets in $\Pi_kM.$

\vspace{5 mm}
\noindent
We begin by considering the equation $\mathcal{R}$ composed by all the elements $X\in\Pi_kM$ whose sources belong to the units of the first groupoid and the targets to the units of the second. Then, of course, any local equivalence will be a local solution of this equation the converse being also true and where a local transformation of the base space operates on $\Pi_kM$ via the conjugation by its $k-$jets. Let us next observe that the first groupoid operates, to the right, on the space $\mathcal{R}$ and that the second groupoid operates to the left. So as to render the data more homogeneous and regular, we shall assume that both these actions are transitive, such an assumptions being much in accordance with the problem considered. In fact, most relevant examples do conform to this requirement and, more important, most contexts where equivalence is sought for also do comply with it. It should be noted that whenever the local equivalence is verified and a local contact transformation is put forward, then one transitivity assumptions gives rise to the other. It now remains to investigate the equation $\mathcal{R}$ as well as its prolongations and determine the desired invariants that will characterize the local equivalence. A first obvious remark is that the dimensions of both groupoids inasmuch as the dimensions of their units sub-spaces must be equal and this also entails the equality, in dimensions, of both the $\alpha$ and the $\beta$ fibres. The standard as well as the merihedric prolongation algorithms do preserve some data and operations but not all. In fact, the prolongation of a groupoid is still, algebraically, a groupoid but its differentiable structure is not guaranteed. Furthermore, the action of the prolonged groupoid on the prolonged equation is not forcibly transitive. On the other hand and concerning the equation $\mathcal{R},$ similar statement can also be affirmed. Nevertheless, it should be emphasized that the above stated transitivity assumption implies that this equation inherits, from either groupoid, a manifold structure rendering it a sub-manifold of $\Pi_kM.$ In what follows, we list several necessary conditions for our problem to have a solution.

\vspace{5 mm}
\noindent
$\mathbf{1.~Differentiability.}$ The iterated prolongations of both groupoids are also Lie groupoids and sub-manifolds of the respective jet spaces.

\vspace{3 mm}
\noindent
$\mathbf{2.~Dimension.}$ Both prolonged groupoids have, at each order, the same dimensions. We do not need to bother about the dimensions of the unit spaces since these spaces are equal (diffeomorphic) to those of the initially given groupoids,

\vspace{3 mm}
\noindent
$\mathbf{3.~Transitivity.}$ Both prolonged groupoids operate, at each order, transitively on the prolonged equation. Consequently, each prolonged equation is a sub-manifold of the corresponding jet space. 

\vspace{3 mm}
\noindent
$\mathbf{4.~Symbols.}$ The symbols of both prolonged groupoids have, at each order, the same dimensions. 

\vspace{3 mm}
\noindent
$\mathbf{5.~\delta-cohomology.}$ Both iterated prolongations of the groupoids have equidimensional $\delta-$cohomologies and thereafter become simultaneously $2-$acyclic (\textit{cf.} \cite{Kumpera1972}).

\vspace{5 mm}
\noindent
We now assume that the above five properties do hold for the given two Lie groupoids of jets. On account of the transitivity of the groupoid actions, we infer that the iterated prolongations of the equation $\mathcal{R}$ become also $2-$acyclic at orders not greater than those for the groupoids. Finally, if we assume further that all the data is real analytic, then Cartan's Theorem\footnote{most improperly called the Cartan-Kähler Theorem} guarantees the desired local equivalence. Unfortunately, we cannot make the same statement in the $C^{\infty}$ realm since, under these weaker assumptions, the whole Cartan Theory breaks down. There are examples of non-analytic though involutive (2-acyclic) equations that do not possess any  solution. In this last realm, much further work is awaiting to be done.

\vspace{5 mm}
\noindent
The local equivalence of analytic Lie groupoids entails the local equivalence of Lie pseudo-groups since these are defined with the help of such groupoids, their defining equations.

\bibliographystyle{plain}
\bibliography{references}

\begin{thebibliography}{10}

\bibitem{Almeida1985}
R.~Almeida and P.~Molino.
\newblock {Suites d'Atiyah et feuilletages transversalement complets}.
\newblock {\em C.R. Acad. Sci. Paris}, 300:13--15, 1985.

\bibitem{Cartan1910}
E.~Cartan.
\newblock {Les systèmes de Pfaff à cinq variables et les équations aux
  dérivées partielles du second ordre}.
\newblock {\em Annales Sci. École Norm. Sup.}, 27:109--192, 1910.

\bibitem{Cartan1914}
E.~Cartan.
\newblock Sur l'équivalence absolue de certains systèmes d'équations
  différentielles et sur certaines familles de courbes.
\newblock {\em Bull. Soc. Math. de France}, 42:12--48, 1914.

\bibitem{Cartan1937}
E.~Cartan.
\newblock Les problèmes d'équivalence.
\newblock {\em Selecta}, pages 113--136, 1937.

\bibitem{Cartan1945}
E.~Cartan.
\newblock {\em Les systèmes différentiels extérieurs et leurs applications
  géométrique}.
\newblock Hermann, Paris, 1945.

\bibitem{Kumpera1999}
A.~Kumpera.
\newblock {On the Lie and Cartan theory of invariant differential systems}.
\newblock {\em J. Math. Sci. Univ. Tokyo}, 6:229--314, 1999.

\bibitem{Kumpera2014}
A.~Kumpera.
\newblock {Automorphisms of Flag Systems}.
\newblock {\em arXiv}, 1411.0951:25, 2015.

\bibitem{Kumpera2017}
A.~Kumpera.
\newblock {On The Equivalence Problem for Geometric Structures, I, II}.
\newblock {\em Journal of Mathematical Sciences, Advances and Applications,
  http://www.scientificadvances.co.in/artical/1/179}, 1412.8394:38, 2015.

\bibitem{Kumpera2016}
A.~Kumpera.
\newblock {Non-integrable Pfaffian systems}.
\newblock {\em arXiv}, 1608.02871:23, 2016.

\bibitem{Kumpera2015}
A.~Kumpera.
\newblock {On the Lie and Cartan theory of invariant differential systems, II}.
\newblock {\em arXiv}, 1511.08703v2:33, 2016.

\bibitem{Kumperaa1982}
A.~Kumpera and C.~Ruiz.
\newblock {Sur l'équivalence locale des systèmes de Pfaff en drapeau}.
\newblock {\em F. Gherardelli (Ed.), Monge-Ampère Equations and Related
  Topics, Firenze 1980, Proceedings, Istit. Naz. Alta Mat. "Francesco Severi",
  Roma}, pages 201--248, 1982.

\bibitem{Kumpera1972}
A.~Kumpera and D.~Spencer.
\newblock {\em Lie Equations, Vol.1: General Theory}.
\newblock Princeton University Press, 1972.

\end{thebibliography}
\end{document}